\documentclass[a4paper]{article}
\usepackage{std-defs,fullrefs,microtype}

\def\epsilon{\varepsilon}
\def\phi{\varphi}

\def\2{{\scriptscriptstyle (2)}} 
\def\P{{\rm P}} 
\def\l{\wedge} 
\def\L{{\cal L}} 
\def\A{{\cal A}} 
\def\S{{\cal S}} 
\def\E{{\cal E}} 

\def\t#1{\tilde{#1}}
\def\wt#1{\widetilde{#1}}
\def\tr#1{{{}^*{#1}\,}}

\def\vecteur#1{{{\vec {#1}}}}

\def\e{{\vecteur e}}
\def\o{{\vecteur o}} 
\def\u{{\vecteur u}}
\def\v{{\vecteur v}}
\def\x{{\vecteur x}}
\def\y{{\vecteur y}}
\def\z{{\vecteur z}}
\def\dz{{\vecteur{\dot z}}}
\def\ddz{{\vecteur{\ddot z}}}

\let\oldsum\sum
\def\sum_#1{\mathop{\oldsum_{#1}}\limits}

\begin{document}

\markboth{Guy Fayolle and Jean-Marc Lasgouttes}
         {A state-dependent polling model with Markovian routing}

\title{A State-Dependent Polling Model \\
       with Markovian Routing}
\author{Guy Fayolle 
       \thanks{Postal address: INRIA --- 
                     Domaine de Voluceau, Rocquencourt, BP105 --- 
                     78153 Le Chesnay, France.}
       \and Jean-Marc Lasgouttes$\ ^*$}
\date{June 1994; revised August 1994, September 1994}
\maketitle

\begin{abstract}
A state-dependent $1$-limited polling model with $N$ queues is
analyzed. The routing strategy generalizes the classical Markovian
polling model, in the sense that two routing matrices are involved,
the choice being made according to the state of the last visited
queue. The stationary distribution of the position of the server is
given. Ergodicity conditions are obtained by means of an associated
dynamical system.  Under rotational symmetry assumptions, average
queue length and mean waiting times are computed.
\end{abstract}

\section{Introduction}

Consider a taxicab in a city in which there are $N$ stations at which
clients arrive and wait for the vehicle. When their turn has come to
be served, they ask for a transit to a destination (one of the other
stations) where they leave the system. Whenever the taxicab finds a
station empty, it goes somewhere else to look for a client. The choice
of the destinations by a client or by the empty taxicab are made via
the two distinct routing matrices $P$ and $\wt P$.

This system can also be seen as a polling model with $N$ queues at
which customers arrive and one server which visits the queues
according to the following rules: if the server finds a client at the
current queue, it serves this client and chooses a new queue according
to the routing matrix $P$; otherwise it selects the next queue
according to $\wt P$. This polling scheme is an extension of the
classical Markovian polling model, with routing probabilities
depending on the state of the last visited queue. As a taxicab only
takes one client at a time, the service strategy is the {\em
$1$-limited\/} strategy, which is known to be more difficult to analyze
than either the {\em gated\/} or {\em exhaustive\/} strategies---in which
the server respectively serves the clients present on its arrival or
serves clients until the queue is empty.

Since the bibliography on polling models is plethoric, we refer the
reader to the references given in Takagi~\cite{Tak:3} and Levy and
Sidi~\cite{LevSid:1}. Among the studies related to our work, we can
cite Kleinrock and Levy~\cite{KleLev:1}, who compute waiting times in
random polling models in which the destination of the server is chosen
independently of its provenance, and Boxma and
Weststrate~\cite{BoxWes:1} who give a pseudo-conservation law for a
model with Markovian routing. Ferguson~\cite{Fer:1} and Bradlow and
Byrd~\cite{BraByr:1} study approximatively a model where the switching
time is station-dependent. Srinivasan~\cite{Sri:1} analyzes a polling
system with the same routing policy as in our model, and with
different service policies at each queue. However, the waiting times
are not computed in the case where the routing is state dependent.
Another important issue concerns state-independent ergodicity
conditions. In this case, necessary and sufficient conditions were
obtained by Fricker and Ja\"\i{}bi~\cite{FriJai:1,FriJai:2} for
deterministic and Markovian routing with various service disciplines.
Borovkov and Schassberger~\cite{BorSch:1} and Fayolle {\em et
al.}~\cite{FayZam:1} give necessary and sufficient conditions for a
system with Markovian routing and 1-limited service.  Finally,
Schassberger~\cite{Sch:2} gives a necessary condition for the
ergodicity of a polling model with 1-limited service and a routing
that depends on the whole state of the system.

Here, we give a new method to get the stationary distribution of the
position of the server and, in the case of a fully symmetrical system,
a way to compute the mean waiting time of a customer. This extends
known results on $1$-limited polling models, especially for the
symmetrical case. However, we do not provide a pseudo-conservation law
as for example Boxma and Weststrate~\cite{BoxWes:1}, since the
computations are not as simple as in the state-independent routing
case. Our method of proof allows to make only minimal assumptions on
the arrival process and applies to either discrete or continuous time
models. Moreover, in the symmetrical case the results allow to compare
various polling strategies.

The paper is organized as follows: in Section~\ref{sec:model}, we
present the model and give a functional equation describing its
evolution. In Section~\ref{sec:proba}, stationary probabilities of the
position of the server are obtained and Section~\ref{sec:ergo} is
devoted to the general problem of ergodicity.  Section~\ref{sec:sym}
presents formulas to compute the first moment of the queue length at
polling instants. These results are used in Section~\ref{sec:wait} to
calculate the waiting time of an arbitrary customer. Applications
to known models are also given.

\section{Description of the model}\label{sec:model}

The system consists of $N$ stations at which clients arrive according
to a stationary process. Let $\S\egaldef\{1,\ldots,N\}$ be the set of
stations. Assume that the server arrives after $n-1$ moves at station
$i\in\S$ where a client is waiting. Then the server loads this client
and goes to station $j\in\S$ with probability $p_{i,j}$, such that
$p_{i,1}+\cdots+p_{i,N}=1$. Conversely, if station $i$ is empty, the
server polls station $j$ with probability $\t p_{i,j}$.  The service
policy is known as {\em 1-limited\/} policy, since at most 1 customer is
served each time a station is visited. The two transition matrices
will be denoted respectively by
\[
P     \;=\; (p_{i,j})_{i,j\in\S},\quad
\wt P \;=\; (\t p_{i,j})_{i,j\in\S}.
\]

The number of new customers arriving to station $q\in\S$ between the
polling of stations $i$ and $j$ when there have been a service (resp.\
no service) at station $i$ is $B_{i,j;q}(n)$ (resp.\ $\wt
B_{i,j;q}(n)$). The vectors
$B_{i,j}(n)=(B_{i,j;1}(n),\ldots,B_{i,j;N}(n))$ (resp.\ $\wt
B_{i,j}(n)$) are i.i.d.\ for different $n$, but their components may
be dependents and not identically distributed.  In the classical
polling terminology, $B_{i,j}(n)$ is the number of arrivals during
{\em both\/} the service and the switchover and $\wt B_{i,j}(n)$ is the
number of arrivals during a switchover.  However, our model covers a
wider class of applications. In particular, the service time may
depend on both $i$ and $j$. In addition, the switchover time
distribution may depend on the state of the last visited queue (empty
or not). Remark also that, unlike in Srinivasan~\cite{Sri:1}, we can
also have $\EE B_{i,j}(n)<\EE\wt B_{i,j}(n)$, which means that serving
a customer can actually take less time than just moving. This would
amount to {\em negative\/} service times in the usual polling
formulation, which are not easy to handle!

It is worth noting that our setting is valid for both discrete-time
and con\-tin\-u\-ous-time models with or without batch arrivals. We do not
need a separate analysis for each case, until the computation of
waiting times, where we will have to give precise definitions of the
vectors $B_{i,j}(n)$ and $\wt B_{i,j}(n)$. We do not describe yet the
exact arrival process and the time taken to travel from one station to
another, since we only need the distribution of the number of clients
arriving during such a move. However, due to the assumption of
independence between successive arrivals, the model applies mainly to
the following situations:
\begin{itemize}
\item discrete-time evolution, when a batch of customers arrives at
each station at the beginning of each time slot; batches are i.i.d.\
with respect to the time slots but they can be correlated between 
stations; 
\item continuous time evolution, when the arrival process is a
compound Poisson process: customers arrive at the system in batches at
the epochs of a Poisson process and the batches have the same
properties as in the previous case.
\end{itemize}

Let $X_i(n)$ be the number of clients waiting at station $i$ at
polling instant $n$ and $S(n)$ be the corresponding position of the
server. If we define $X(n)\egaldef(X_1(n),\ldots,X_N(n))$, then
$\L\egaldef(S,X)$ forms a Markov chain. Throughout this paper, we will assume
that this chain is irreducible. For any vector
$\z=(z_1,\ldots,z_N)\in{\cal D}^N$ (where ${\cal D}$ denotes the unit
disc in the complex plane), we note $\z^B=z_1^{B_1}\cdot
z_2^{B_2}\cdots z_N^{B_N}$ and define the following generating
functions:
\begin{eqnarray*}
a_{i,j}(\z)    & \egaldef & \EE[\z^{B_{i,j}(n)}],\\
\t a_{i,j}(\z) & \egaldef & \EE[\z^{\wt B_{i,j}(n)}],\\
F_i(\z;n)      & \egaldef & \EE[\z^{X(n)}\ind{S(n)=i}],\\
\wt F_i(\z;n)  & \egaldef & \EE[\z^{X(n)}\ind{S(n)=i,X_i(n)=0}],\\
      & = & F_i(z_1,\ldots,z_{i-1},0,z_{i+1},\ldots,z_N;n),\\
F_i(\z)        & \egaldef & \lim_{n\to\infty}F_i(\z;n),\\
\wt F_i(\z)    & \egaldef & \lim_{n\to\infty}\wt F_i(\z;n),
\end{eqnarray*}
where $\1_{\cal E}$ is as usual the indicator function of the set
${\cal E}$.  The following result holds:

\begin{thm}
Let $A(\z)$, $\wt A(\z)$ and $\Delta(\z)$ be matrices defined as
follows: for all $i,j\in\S$, the elements of row $i$ and column $j$
are given by
\begin{eqnarray*}
[A(\z)]_{i,j}&=&p_{j,i}a_{j,i}(\z),\\
{}[\wt A(\z)]_{i,j}&=&\t p_{j,i}\t a_{j,i}(\z),\\
{}[\Delta(\z)]_{i,j}&=& {1\over z_i}\ind{i=j}.
\end{eqnarray*}

\noindent Then the vectors
$F(\z)=(F_1(\z),\ldots,F_N(\z))$ and $\wt F(\z)=(\wt
F_1(\z),\ldots,\wt F_N(\z))$ are related by the functional equation
\begin{equation}\label{eq:matrix}
[I-A\Delta(\z)]F(\z)=[\wt A(\z)-A\Delta(\z)]\wt F(\z),
\end{equation}
where $A\Delta(\z)$ stands for $A(\z)\Delta(\z)$. 
\end{thm}

\begin{proof}{}
We have, for all $i,j\in\S$ and $n>0$,
\begin{eqnarray*}
\lefteqn{ \EE[\z^{X(n+1)}\ind{S(n)=i,S(n+1)=j}]}\qqqq\\
&=&  \EE[\z^{X(n)+B_{i,j}(n)-\e_i}\ind{S(n)=i,S(n+1)=j,X_i(n)>0}]\\
& & +\EE[\z^{X(n)+\wt B_{i,j}(n)}\ind{S(n)=i,S(n+1)=j,X_i(n)=0}]\\
\smallskip
&=& p_{i,j}a_{i,j}(\z)\EE[\z^{X(n)-\e_i}\ind{S(n)=i,X_i(n)>0}]\\
& & +\t p_{i,j}\t a_{i,j}(\z)\EE[\z^{X(n)}\ind{S(n)=i,X_i(n)=0}]\\
\smallskip
&=& p_{i,j}a_{i,j}(\z){F_i(\z;n)-\wt F_i(\z;n)\over z_i}
    +\t p_{i,j}\t a_{i,j}(\z)\wt F_i(\z;n).
\end{eqnarray*}

The second equality above uses the independence of the routing and the 
arrivals with respect to the past. Summing over all possible $i$, we get
\[
  F_j(\z;n+1)
  = \sum_{i=1}^N p_{i,j}a_{i,j}(\z){F_i(\z;n)-\wt F_i(\z;n)\over z_i}
    +\sum_{i=1}^N \t p_{i,j}\t a_{i,j}(\z)\wt F_i(\z;n),
\]
so that, letting $n\to\infty$,
\[
  F_j(\z)-\sum_{i=1}^N{p_{i,j}a_{i,j}(\z)\over z_i}F_i(\z)
  = \sum_{i=1}^N\t p_{i,j}\t a_{i,j}(\z)\wt F_i(\z)
   -\sum_{i=1}^N{p_{i,j}a_{i,j}(\z)\over z_i}\wt F_i(\z),
\]

which is equivalent to~(\ref{eq:matrix}).
\end{proof}

Defining, when it exists,
\begin{equation}\label{eq:defD}
 D(\z)\egaldef[I-A\Delta(\z)]^{-1}[\wt A(\z)-A\Delta(\z)],
\end{equation}
one sees that~(\ref{eq:matrix}) can be rewritten as
\begin{equation}\label{eq:matrix2}
F(\z)=D(\z)\wt F(\z).
\end{equation}

This functional equation contains all the information sufficient to
characterize $F(\z)$ and $\wt F(\z)$. Although its solution for $N\geq
3$ is still an open question, partial results can be derived as
shown in the following sections.

\section{The stationary distribution of the position of the server}
\label{sec:proba}

What renders the polling model presented in the previous section
difficult to analyze is, among other things, the fact that the
movements of the server depend on the state of the visited
stations. In particular, $\{S(n)\}_{n\geq0}$ is not a Markov process,
as it would be when $p_{i,j}=\t p_{i,j}$ for all $i$ and $j$. Some
computations are needed to get the stationary probabilities of the
position of the server, that is
\begin{eqnarray*}
F_i(\e)     &=& P(S=i) \\
\wt F_i(\e) &=& P(S=i,X_i=0).
\end{eqnarray*}

It appears that these stationary probabilities depend not only on the
transition probabilities, but also on the following mean values:
\begin{eqnarray}
\alpha_{i;q} 
  &\egaldef& \sum_{j=1}^N p_{i,j}\EE B_{i,j;q}(n),
                                                  \label{eq:defaliq}\\
\t\alpha_{i;q} 
  &\egaldef& \sum_{j=1}^N \t p_{i,j}\EE\wt B_{i,j;q}(n).
                                                  \label{eq:deftaliq}
\end{eqnarray}

Here $\alpha_{i;q}$ (resp.\ $\t\alpha_{i;q}$) is the mean number of
new clients at station $q$ between the arrival of the server at
station $i$ which was non-empty (resp.\ empty) and its arrival at the
next (arbitrary) polled station.  Since the arrival process is
stationary, it is possible to define for all $q\in\S$ the mean number
$\lambda_q$ of arrivals at station $q$ per unit of time and write,
with obvious definitions for $\tau_i$ and $\t\tau_i$,
\begin{eqnarray}
 \alpha_{i;q} &=& \lambda_q\tau_i, \label{eq:defaliq2}\\
 \t\alpha_{i;q} &=& \lambda_q\t\tau_i. \label{eq:deftaliq2}
\end{eqnarray}

Throughout this paper we will sometimes write $F$ (resp.\ $\wt F$)
instead of $F(\e)$ (resp.\ $\wt F(\e)$) in order to shorten the
notation. Moreover, $\tr F$ denotes the transpose of the vector $F$.

\begin{thm}\label{thm:proba}
Define the matrices $\A\egaldef (\alpha_{i;q})$ and $\wt\A\egaldef
(\t\alpha_{i;q})$.  When the
Markov chain $\L$ is ergodic, $F$ and $\wt F$ satisfy the following
linear system of equations:
\begin{eqnarray}
\tr F\e 
  &=& 1,\label{eq:FFtnormal}\\
\tr F[I-P] 
  &=& \tr{\wt F}[\wt P-P],\label{eq:FFtproba}\\
\tr F[I-\A] 
  &=& \tr{\wt F}[I-\A+\wt\A].\label{eq:FFtflux}
\end{eqnarray}

Define
\[
\hat\rho\egaldef\sum_{i=1}^N\lambda_i(\tau_i-\t\tau_i).
\]

When $\hat\rho\neq 1$, Equations~(\ref{eq:FFtproba})
and~(\ref{eq:FFtflux}) can be rewritten as
\begin{eqnarray}
F_j &=&
  \sum_{i\in\S}\t p_{i,j}F_i+
  \bar\tau\sum_{i\in\S}\lambda_i(p_{i,j}-\t p_{i,j})
  \mbox{,\ \ \ for all }j\in\S,\label{eq:resF}\\
  \bar\tau &=& {1\over 1-\hat\rho}\sum_{j\in\S}F_j\t\tau_j,\label{eq:restau}
\end{eqnarray}
and $\wt F$ is given by
\begin{equation}\label{eq:resFt}
\wt F_j = F_j - \lambda_j\bar\tau.
\end{equation}

Moreover the mean time between two consecutive visits to queue $j$ is 
\begin{equation}\label{eq:cycle}
\EE T_{j,j} = {\bar\tau\over F_j}
\end{equation}
\end{thm}

The quantity $\bar\tau$ defined in~(\ref{eq:restau}) can be seen as
the mean time between two polling instants. Note that, with the above
definitions, $\hat\rho$ is in general different from the classical
$\rho$ defined in queueing theory. However, they coincide if $P=\wt P$
and the time between two polling instants can be decomposed into a
state-independent switchover time and a service time which only
depends on the station where the customer is.

\begin{proof}{ of Theorem~\ref{thm:proba}}
Although this theorem could be proved analytically using the
functional equation~(\ref{eq:matrix}), we give here a simple
probabilistic interpretation of~(\ref{eq:FFtproba})
and~(\ref{eq:FFtflux}). Indeed, when the system is
ergodic,~(\ref{eq:FFtproba}) follows directly from
\[
P(S=j) = \sum_{i=1}^N \Bigl[P(S=i,X_i>0)p_{i,j}+P(S=i,X_i=0)\t p_{i,j}\Bigr].
\]

Moreover, writing the equality of the outgoing and ingoing flows at
station $q$ gives
\begin{eqnarray*}
P(S=q,X_q>0) &=& \sum_{i=1}^N\sum_{j=1}^N 
  \Bigl[P(S=i,X_i>0)p_{i,j}\EE B_{i,j;q}(n)\\
& &   \mbox{}+P(S=i,X_i=0)\t p_{i,j}\EE \wt B_{i,j;q}\Bigr],
\end{eqnarray*}
which is equivalent to~(\ref{eq:FFtflux}) if we take into
account~(\ref{eq:defaliq}) and~(\ref{eq:deftaliq}).  
It is straightforward to see, using~(\ref{eq:defaliq2})
and~(\ref{eq:deftaliq2}), that, when $\hat\rho\neq1$,
\[
 [I-\A][I-\A+\wt\A]^{-1}=I-{\wt\A\over 1-\hat\rho}.
\]

With this relation and $\bar\tau$ as in~(\ref{eq:restau}),
Equations~(\ref{eq:resF}) and~(\ref{eq:resFt}) follow easily from
(\ref{eq:FFtproba}) and~(\ref{eq:FFtflux}). Let $N_{j,j;j}$ be the
number of arrivals to a queue $j$ between two consecutive visits of
the server.  The equality of flows reads
\[
 P(X_j>0\mid S=j) = \EE N_{j,j;j}.
\]

Since the last equation can be rewritten as
\[
 {F_j-\wt F_j\over F_j}=\lambda_j\EE T_{j,j}\ ,
\]
we obtain~(\ref{eq:cycle}) and the proof of the theorem is concluded.
\end{proof}

The set of
equations~(\ref{eq:FFtnormal}),~(\ref{eq:resF}),~(\ref{eq:restau})
and~(\ref{eq:resFt}) provides a convenient way to compute $F$ and $\wt
F$.  The following proposition gives simple conditions under which its
rank is full.

\begin{pro}\label{pro:rank}
The linear system given by~(\ref{eq:FFtnormal}),~(\ref{eq:resF})
and~(\ref{eq:restau}) has a unique solution when $\wt P$ has exactly
one essential class. When $\wt P$ has $K>1$ essential classes
$\E_1,\ldots,\E_K$, the Markov chain $\L$ is never ergodic unless
\begin{equation}\label{eq:compat}
 \sum_{j\in\E_m}\sum_{i\in\S}\lambda_i(p_{i,j}-\t p_{i,j})=0 
  \mbox{,\ \ for all }m\leq K.
\end{equation}
\end{pro}

\begin{proof}{}
Define $\psi_i\egaldef F_i/\bar\tau$. Then the
system~(\ref{eq:resF})--(\ref{eq:restau}) reads
\begin{eqnarray}
 \psi_j-\sum_{i\in\S}\t p_{i,j}\psi_i
  &=& \sum_{i\in\S}\lambda_i(p_{i,j}-\t p_{i,j})
      \mbox{,\ \ \ for all }j\in\S,\label{eq:psi1}\\
 \sum_{i\in\S}\psi_i\t\tau_i &=& 1-\hat\rho\label{psi2}.
\end{eqnarray}

The system~(\ref{eq:psi1}) have rank $N-1$ unless the stochastic
matrix $\wt P$ has several essential classes.  Then, for each $m\leq
K$, $i\in\E_m$ and $j\not\in \E_m$, $\t p_{i,j}=0$. In this
case,~(\ref{eq:psi1}) has solutions only when~(\ref{eq:compat}) holds.
\end{proof}

The conditions~(\ref{eq:compat}) represent in some sense ``zero
drift'' relationships which, as usual, imply more involved
derivations.  However, when the arrivals form independent compound
Poisson processes at each queue, the system is never ergodic!  The
proof is of analytic nature. It requires Taylor expansions of second
order in equation~(\ref{eq:matrix}), similar to those extensively used
in Section~\ref{sec:sym}. In general, when the batches are correlated,
it is difficult to conclude and we suspect that all situations might
occur (ergodicity, null recurrence or transience).

It is worth noting the special role played by $\wt P$. In particular,
when the arrivals are Poisson, the Markov chain is {\em never\/}
ergodic if $\wt P$ admits several essential classes (obviously, $P$
has then to be chosen to ensure the irreducibility of $\L$). In our
opinion, this result is not intuitive and we cannot be explained just
by waving hands.

\section{Conditions for ergodicity}
\label{sec:ergo}

The purpose of this section is to classify the process $\L$, viewed as
a random walk on $\S\times\ZZ^N_+$, in terms of ergodicity and
non-ergodicity.

\subsection{Necessary condition}\label{sub:nec}

We give a necessary condition for ergodicity which is a simple
consequence of the results of Theorem~\ref{thm:proba}. When the system
is ergodic, we have $\wt F_i>0$ for all $i\in\S$ and
Equation~(\ref{eq:resFt}) implies the following.

\begin{thm}\label{thm:nec}      
If the Markov chain $\L$ is ergodic, then
\begin{equation}\label{eq:nec:rho}
  \hat\rho<1,
\end{equation}
\begin{equation}\label{eq:nec:flux}
  \lambda_i\bar\tau<F_i\mbox{\rm ,\ \ for all }i\in\S.
\end{equation}
\end{thm}

These conditions extend in the state-dependent case the results
obtained in~\cite{BorSch:1}.  When they hold, we prove another useful
inequality.  Indeed, instantiating~(\ref{eq:nec:flux})
in~(\ref{eq:restau}) yields
\[
 \bar\tau> {\sum_{j\in\S}\bar\tau\lambda_j\t\tau_j\over 1-\hat\rho},
\]
or, equivalently,
\begin{equation}\label{eq:nec:sup}
  1-\sum_{j\in\S}\lambda_j\tau_j> 0.
\end{equation}

\subsection{Sufficient condition}

Let us emphasize that this part could be skipped by readers not
strongly interested (if any!) in ergodicity. The mathematical
understanding requires the reading of Appendix~\ref{app:induc}, which
summarizes deep works of~\cite{MalMen:1},~\cite{FayMalMen:1}
and~\cite{FayZam:1}. The approach relies on the study of a
dynamical system associated to $\L$. For the sake of readability, we
recall hereafter two main notions.

\begin{itemize}
\item {\em Induced chain $\L^\l$}. It is a Markov chain
corresponding to a polling system in which the queues belonging to the
{\em face\/} $\l\subset\S$ are kept saturated. From a purely
notational point of view, the original system would correspond to the
case $\l=\emptyset$. The behaviours of $\L$ and $\L^\l$ are not {\em
directly\/} connected to each other.
\item {\em Second vector field $\v^\l$}. One can imagine that the random
walk starts from a point which is close to $\l$, but sufficiently far
from all other faces $\l'$, with $\l\not\subset\l'$. After some
time---sufficiently long, but less than the minimal distance from
$\l'$---, the stationary regime in the induced chain will be
installed. In this regime, one can ask about the mean along $\l$: it
is defined exactly by $\v^\l$.
\end{itemize}

As in Theorem~\ref{thm:proba}, the stationary position of the server
for any ergodic induced chain can be obtained as the solution of a
linear system consisting of $N+1$ equations.  Indeed, for all
$t\in\S$, $j\not\in\l$, we have
\begin{eqnarray}
\pi^\l(t)
  &=& \sum_{s\in\l}\pi^\l(s)p_{s,t}
     +\sum_{s\not\in\l}\pi^\l(s)\t p_{s,t}\label{eq:prob1}\\
  & &  +\mbox{}\sum_{s\not\in\l}\pi^\l(s,x_s>0)(p_{s,t}-\t p_{s,t}),
                                                                 \nonumber\\
\pi^\l(j,x_j>0) 
  &=& \sum_{s\in\l}\pi^\l(s)\lambda_j\tau_s
     +\sum_{s\not\in\l}\pi^\l(s)\lambda_j\t\tau_s\label{eq:flux1}\\
  & &\mbox{}+\sum_{s\not\in\l}\pi^\l(s,x_s>0)\lambda_j(\tau_s-\t\tau_s),
                                                                  \nonumber\\
\sum_{s\in\S}\pi^\l(s) &=& 1.\label{eq:norm1}
\end{eqnarray}

Let us introduce the quantities
\begin{eqnarray*}
\hat\rho^\l 
  &\egaldef& \sum_{s\not\in\l}\lambda_s(\tau_s-\t\tau_s),\\
\bar\tau^\l &\egaldef&
   {1\over 1-\hat\rho^\l}
       \biggl[\sum_{s\in\l}\pi^\l(s)\tau_s
       +\sum_{s\not\in\l}\pi^\l(s)\t\tau_s\biggr].
\end{eqnarray*}

Note that when $\l$ is ergodic,
$\pi^\l(j,x_j>0)=\lambda_j\bar\tau^\l$, for all $j\not\in\l$, and
$\bar\tau^\l$ can be interpreted as the mean time between two polling
instants of the induced chain.  Then the system defined
by~(\ref{eq:prob1}) and~(\ref{eq:flux1}) can be replaced by the
forthcoming $N+1$ equations: for all $t\in\S$,
\begin{eqnarray}
\pi^\l(t)
  &=& \sum_{s\in\l}\pi^\l(s)p_{s,t}
       +\sum_{s\not\in\l}\pi^\l(s)\t p_{s,t}
       +\bar\tau^\l\sum_{s\not\in\l}\lambda_s(p_{s,t}-\t p_{s,t}),
                                                        \label{eq:prob2}\\
\bar\tau^\l
  &=& {1\over 1-\hat\rho^\l}
       \biggl[\sum_{s\in\l}\pi^\l(s)\tau_s
       +\sum_{s\not\in\l}\pi^\l(s)\t\tau_s\biggr].
       \label{eq:flux2}
\end{eqnarray}

When $\l=\emptyset$, the
system~(\ref{eq:norm1}),~(\ref{eq:prob2}) and~(\ref{eq:flux2})
coincides with~(\ref{eq:FFtnormal}),~(\ref{eq:resF})
and~(\ref{eq:restau}). In addition, under~(\ref{eq:nec:sup}), we have
\begin{equation}\label{eq:rhol}
\hat\rho^\l<1,\mbox{ for all }\l\mbox{'s}.
\end{equation}

An easy flow computation shows that the components of the drift vector
$\vecteur M(s,\x)$ can be expressed as
\[
M_j(s,\x)\;=\; \lambda_j\tau_s\ind{x_s>0}
                     +\lambda_j\t\tau_s\ind{x_s=0}
                     -\ind{x_s>0,s=j}.
\]

Then the computation of the second vector field becomes easy.  For
$j\in\l$,
\begin{eqnarray*}
v^\l_j &=& \sum_{x\in C^\l}_{s\in \S} \pi^\l(s,\x)M_j(s,\x)\\
  &=& \sum_{x\in C^\l}_{s\not\in\l} \pi^\l(s,\x)M_j(s,\x)
     +\sum_{x\in C^\l}_{s\in\l} \pi^\l(s,\x)M_j(s,\x)\\
  &=& \lambda_j\sum_{s\not\in\l}\Bigl[\pi^\l(s,x_s>0)\tau_s
                            +\pi^\l(s,x_s=0)\t\tau_s\Bigr]\\
  & &\mbox{}+\lambda_j\sum_{s\in\l}\pi^\l(s)\tau_s-\pi^\l(j)\\
  &=&\lambda_j\bar\tau^\l-\pi^\l(j).
\end{eqnarray*}

In order to apply Theorem~\ref{thm:FayZam}, it will be convenient
to introduce
\begin{eqnarray}
f_i(\x) 
  &\egaldef& \Bigl[\tr\x[I-\A+\wt \A]^{-1}\Bigr]_i\nonumber\\
  &=&        x_i+{\lambda_i\sum_{j=1}^Nx_j(\tau_j-\t\tau_j)
                   \over1-\hat\rho}.\label{eq:lyapi}
\end{eqnarray}

This function is not as outlandish as it could seem at first sight. It
is indeed directly related to flow conservation equations
(see~\cite{FayZam:1}).  Let $\l$ be any ergodic face. For $i\in\l$,
we get after a straightforward computation
\begin{eqnarray*}
f_i(\v^\l)
  &=& v^\l_i
      +{\lambda_i\over1-\hat\rho}\sum_{j=1}^N(\tau_j-\t\tau_j)v^\l_j\\
  &=& {\lambda_i\sum_{j=1}^N\pi^\l(j)\t\tau_j\over1-\hat\rho}
      -\pi^\l(i).
\end{eqnarray*}

For $i\not\in\l$, and from the very definition of $f_i(\x)$, we have
\[
f_i(\v^\l)-\sum_{j=1}^Nf_j(\v^\l)\lambda_i(\tau_j-\t\tau_j)=v^\l_i=0,
\]
or
\[
 f_i(\v^\l)={\lambda_i\sum_{j\in\l}f_j(\v^\l)(\tau_j-\t\tau_j)
             \over1-\hat\rho^\l}.
\]

\medskip
\begin{thm}\label{thm:suff}
Assume that~(\ref{eq:nec:rho})--(\ref{eq:nec:flux}) hold and that, for
any ergodic face $\l$,
\begin{equation}\label{eq:suf}
 f_i(\v^\l)\equiv
  {\lambda_i\sum_{j=1}^N\pi^\l(j)\t\tau_j\over1-\hat\rho}
      -\pi^\l(i) <0,\ \ \mbox{for all }i\in\l.
\end{equation}

Then the random walk $\L$ is ergodic.  In particular, when $P=\wt P$,
the conditions~(\ref{eq:nec:rho}) and~(\ref{eq:nec:flux}) are
necessary and sufficient for the random walk $\L$ to be ergodic.
\end{thm}

\begin{proof}{}
We shall apply Theorem~\ref{thm:FayZam} (quoted in
Appendix~\ref{app:induc}), which in itself contains a principle of
gluing Lyapounov functions together. Most of the time, these
functions are piecewise linear. Since the $f_i(\v^\l)$ enjoy ``nice''
properties, it is not unnatural to search for a linear combination
like
\[
  f(\x)\egaldef \sum_{i=1}^N u_if_i(\x),
\]
where $\u=(u_1,\ldots,u_N)$ is a positive vector to be properly determined.
Then, for any ergodic $\l$, 
\[
  f(\v^\l)=\sum_{i\in\l}\Bigl[u_i+
   {\sum_{j\not\in\l}\lambda_ju_j\over1-\hat\rho^\l}(\tau_i-\t\tau_i)\Bigr]
   f_i(\v^\l).
\]

The basic constraint for $\u$ is to ensure the positivity of $f$.
Using~(\ref{eq:nec:sup}) and~(\ref{eq:rhol}), it appears that a
suitable choice is $u_i=\max(\t\tau_i-\tau_i,\epsilon)$, for some
$\epsilon$, positive and sufficiently small.  Then one can directly
check that $f(\x)>0$, for any $\x\in\RR^N_+$, and $f(\v^\l)<0$, for
any ergodic face $\l$.
\end{proof}

Although the conditions of Theorem~\ref{thm:suff} are sufficient for
ergodicity, they might well be implied by
conditions~(\ref{eq:nec:rho})--(\ref{eq:nec:flux}) only. We did not
check this fact, since it is difficult to compare algebraically the
solutions of system~(\ref{eq:prob2})--(\ref{eq:flux2}), for different
$\l$'s. Let us simply formulate the conjecture
that~(\ref{eq:nec:rho})--(\ref{eq:nec:flux}) are also sufficient for
ergodicity when $P\neq\wt P$.

\begin{thm}
Assume that $P=\wt P$ and that one of the following hold:
\begin{enumerate}
\item there exists $j\in\S$, such that
$\lambda_j\bar\tau-F_j>0$;
\item $\hat\rho>1$.  
\end{enumerate}

Then the random walk $\L$ is transient.
\end{thm}

\begin{proof}{}
We assume that the queues are numbered according to the following
order:
\begin{equation}\label{eq:order}
{\lambda_1\over F_1}\leq{\lambda_2\over F_2}\leq\ldots\leq
{\lambda_N\over F_N}.
\end{equation}

Consider all faces $\l_i\egaldef\{i,\ldots,N\}$, for
$i=1,\ldots,N$. The idea is to show the transience of $\L$ by visiting
the ordered set $\l_1,\l_2,\ldots\l_N$. In the usual terminology of
dynamical systems, this amounts to finding a set of trajectories going
to infinity with positive probability. To this end, the following
algebraic relationship will be useful: for any $\l$ and $k\not\in\l$,
we have
\begin{equation}\label{eq:compare}
(\lambda_k\bar\tau^{\l+\{k\}}-F_k)(1-\hat\rho^{\l+\{k\}})
 = (\lambda_k\bar\tau^\l-F_k)(1-\hat\rho^\l).
\end{equation}

By definition, the face $\l_1\equiv\S$ is ergodic (see
Appendix~\ref{app:induc}). Assume now $\l_i$ is ergodic, for some
fixed $i\in\S$. If $v^{\l_i}_i<0$, then one can prove the inequality
\[
 1-\sum_{s=1}^i\lambda_s\tau_s
   \leq {1\over\bar\tau^{\l_i}}\sum_{s=i}^NF_s\tau_s,
\]
which in turn yields $1-\hat\rho^{\l_i}>0$. By
using~(\ref{eq:compare}),~(\ref{eq:order}) and Theorem~\ref{thm:suff},
it follows that the face $\l_{i+1}$ is ergodic.  Conversely, if
$v^{\l_i}_i>0$, then $\v^{\l_i}>\vecteur0$ and the random walk $\L$ is
transient (see~\cite{FayMalMen:1}).

Thus, by induction, we have shown that either $\L$ is transient or the
face $\l_N\equiv\{N\}$ is ergodic. In the latter case, the assumptions
of the theorem, together with~(\ref{eq:compare}) and~(\ref{eq:order}),
yield $v^{\l_N}_N\equiv\lambda_N\bar\tau^{\l_N}-F_N>0$, so that
$\L$ is transient.
\end{proof}

\subsection{Remarks and problems}

\begin{enumerate}
\item In general, the second vector field requires to compute 
invariant measures of Markov chains in dimensions strictly smaller
than $N$ (the induced chains). This is rarely possible for $N>3$, but
there are some miracles, the most noticeable ones concerning Jackson
networks and some polling systems~\cite{FayZam:1}. 
\item No stochastic monotonicity argument seems to work, when the
$\tau_k-\t\tau_k$'s are not all of the same sign, even for $P=\wt
P$. This monotonicity exists and is crucial
in~\cite{BorSch:1},~\cite{FriJai:1} and~\cite{FriJai:2}. 
\item The analysis of the vector field $\v^\l$ is interesting in
itself and will be the subject of a future work. Let us just quote one
negative property: when $\l$ is ergodic and $\v_k^\l<0$, then
$\l_1=\l\setminus\{k\}$ can be non-ergodic, contrary to what happens
in Jackson networks (see~\cite{FayZam:1}). In the case $P=\wt P$,
we conjecture nonetheless that the dynamical system formed from the
velocity vectors $\v^\l$ is strongly acyclic (in the sense that the same
face is not visited twice).
\end{enumerate}

\section{The first moment of the queue length}\label{sec:sym}

This section shows that, when the system enjoys some symmetry
properties, it becomes possible to derive the stationary mean queue
length $\EE [X_m\mid S=m]$ seen by the server at polling instants.  To
this end, ergodicity of $\L$ will be assumed, as well as the existence
of second moments of all random variables of interest.

\begin{ass}\label{ass:sym}
The system has a {\em rotational\/} symmetry
\begin{enumerate}
\item for all $i,j\in\S$, 
  $p_{i,j}=p_{j-i}$, $\t p_{i,j}=\t p_{j-i}$, 
  $a_{i,j}(\z)=a_{j-i}(\z)$ and $\t a_{i,j}(\z)=\t a_{j-i}(\z)$;
\item for all $1\leq s,d\leq N$, 
  \begin{eqnarray*}
    a_d(z_1,\ldots,z_N)    &=& a_d(z_{1+s},\ldots,z_{N+s}),\\
    \t a_d(z_1,\ldots,z_N) &=& \t a_d(z_{1+s},\ldots,z_{N+s}).
  \end{eqnarray*}
\end{enumerate}
\end{ass}

Note that this assumption is less restrictive than the ones appearing
in the literature, even for state-independent situations.  It allows
in particular cyclic and random polling
strategies. From~\ref{ass:sym}, we directly get, for all $1\leq
i,s\leq N$,
\begin{eqnarray}
F_i(\z)
 &=&  F_{i+s}(z_{1+s},\ldots,z_{N+s}),\label{eq:symF}\\
\wt F_i(\z)
 &=&\wt F_{i+s}(z_{1+s},\ldots,z_{N+s}).\label{eq:symFt}
\end{eqnarray}

Consequently, $P(S=i)=F_i(\e)=1/N$. In addition, $A(\z)$ and $\wt
A(\z)$ are {\em circulant\/} matrices---see Appendix~\ref{app:circu}
for the results and notation which will be used from now on. It
implies in particular that, for $k\in\S$, $\v_k$ is an eigenvector of
$A(\z)$ and $\wt A(\z)$ with respective eigenvalues
\begin{eqnarray*}
\mu_k(\z)   &\egaldef& \sum_{d=1}^N p_da_d(\z)\omega_k^d\,,\\
\t\mu_k(\z) &\egaldef& \sum_{d=1}^N \t p_d\t a_d(\z)\omega_k^d.
\end{eqnarray*}

Note that $\mu_N(\z)$ (resp.\ $\t\mu_N(\z)$) is the generating
function of the number of arrivals during the transportation of one
client (resp.\ a move without client). Moreover,
$\mu_k\egaldef\mu_k(\e)$ and $\t\mu_k\egaldef\t\mu_k(\e)$ are
eigenvalues of $P$ and $\wt P$. For {\em a priori\/} technical reasons
we shall also need

\begin{ass}\label{ass:irred}
For any $k<N$, $\t\mu_k\neq1$.
\end{ass}

These relations are not surprising in view of
Proposition~\ref{pro:rank}. Since the system is symmetrical, the mean
number of customers arriving during a travel of the taxicab
$\alpha_{i;q}$ and $\t\alpha_{i;q}$ as defined
by~(\ref{eq:defaliq})--(\ref{eq:deftaliq}) does not depend on $i$ and
$q\,$; we note them $\alpha$ and $\t\alpha$ and remark that
\begin{eqnarray*}
 \alpha   &=& \d{z_m}{\mu_N}(\e),\\
 \t\alpha &=& \d{z_m}{\t\mu_N}(\e).
\end{eqnarray*}

The second derivatives of $\mu_N(\z)$ and $\t\mu_N(\z)$ are defined as
\begin{eqnarray}
 \alpha^\2_{q,r}&\egaldef&
  {\partial^2\mu_N(\z)\over\partial z_q\partial z_r}(\e),\label{eq:al2}\\
 \t\alpha^\2_{q,r}&\egaldef&
  {\partial^2\t\mu_N(\z)\over\partial z_q\partial z_r}(\e)\label{eq:tal2}
\end{eqnarray}

Note that these values depend only of the unsigned distance between
$q$ and $r$ and that we can write
$\alpha^\2_{q,q+d}=\alpha^\2_{q+d,q}\egaldef\alpha^\2_d$,
$\t\alpha^\2_{q,q+d}=\t\alpha^\2_{q+d,q}\egaldef\t\alpha^\2_d$.  It
will be also convenient to introduce the first and second moments of
the number of arrivals between two polling instants, that is
\begin{eqnarray}
\bar\alpha 
  &\egaldef& P(X_m>0\mid S=m)\alpha+\P(X_m=0\mid S=m)\t\alpha,\label{eq:bal}\\
\bar\alpha^\2_d 
  &\egaldef& P(X_m>0\mid S=m)\alpha^\2_d
             +\P(X_m=0\mid S=m)\t\alpha^\2_d\label{eq:bal2}.
\end{eqnarray}

Using the symmetry of the system and according to the properties of
generating functions, we have for any $m\in\S$
\begin{equation}\label{eq:mean}
 \EE [X_m\mid S=m]=N\d{z_m}{F_m}(\e).
\end{equation}

By symmetry, this expression is independent from $m$. In order to
derive $\EE[X_m\mid S=m]$, we use the fact that the matrices
$I-A\Delta(\z)$ and $\wt A(\z)-A\Delta(\z)$ are not invertible when
$\z=\e$. This implies that the matrix $D(\z)$ as defined
in~(\ref{eq:defD}) is not continuous in a neighborhood of $\z=\e$ and
provides a mean to compute the derivatives of $F(\z)$ and $\wt F(\z)$
for $\z=\e$. However in our setting, even the existence of $D(\z)$ in
a neighborhood of $\e$ is difficult to prove and the computations
become really involved. The approach that we present here avoids the
theoretical problems and simplifies the computations. In the next
lemma, we show how it is possible to choose $\z$ such that the left
member of Equation~(\ref{eq:matrix}) can be rendered of ``small''
order w.r.t.\ the right one.

\begin{lem}\label{lem:ut}
Let $t\to\z_t$ be a function from $\RR_-\to{\cal D}^N$ such that in a
neighborhood of $t=0$, we can write
$\z_t=\e+t\dz+{1\over2}t^2\ddz+\o(t^2)$, where 
$\dz=\dot z_1\e_1+\cdots+\dot z_N\e_N
=\dot\zeta_1\v_1+\cdots+\dot\zeta_N\v_N\in\CC^N$ 
and $\ddz=(\ddot z_1,\ldots,\ddot z_N)\in\CC^N$. 
Assume that $\dot z_1+\cdots+\dot z_N=0$. Then there exists for $t>0$
a vector $\u_t$ such that, for $1\leq k<N$,
\begin{eqnarray}
\tr{\u_t}[I-A\Delta(\z_t)]\v_k
 &=& o(t),\label{eq:dltuvk}\\
\tr{\u_t}[\wt A(\z_t)-A\Delta(\z_t)]\v_k
 &=& t{1-\t\mu_k\over1-\mu_k}
  \dot\zeta_{N-k}+o(t),\label{eq:dltuAtvk}
\end{eqnarray}
and
\begin{eqnarray}
\tr{\u_t}[I-A\Delta(\z_t)]\e
 &=& -t^2\sum_{l=1}^{N-1}\dot\zeta_l\dot\zeta_{N-l}
  \biggl[{1\over 1-\mu_l}
  +{N\over2}\sum_{d=1}^N\alpha^\2_d\omega_l^d\biggr]\nonumber\\
 & &\mbox{} +t^2{1-N\alpha\over2N}\sum_{q=1}^N\ddot z_q+o(t^2),
                                                        \label{eq:dltue}\\
\tr{\u_t}[\wt A(\z_t)-A\Delta(\z_t)]\e
 &=&  -t^2\sum_{l=1}^{N-1}\dot\zeta_l\dot\zeta_{N-l}
    \biggl[{1\over 1-\mu_l}
    +{N\over2}\sum_{d=1}^N(\alpha^\2_d-\t\alpha^\2_d)\omega_l^d\biggr]
                                                        \nonumber\\
 & &\mbox{}+t^2{1-N\alpha+N\t\alpha\over2N}\sum_{q=1}^N\ddot z_q+o(t^2).
                                                        \label{eq:dltuAte}
\end{eqnarray}
\end{lem}

\begin{proof}{}
See Appendix~\ref{app:proof}.
\end{proof}

For the sake of simplicity, the computations have been carried out in
the most natural way which apparently requires $\mu_k\neq1$ for $k<N$.
In fact, the reader can convince himself that all derivations could be
achieved by rendering the right member of~(\ref{eq:matrix}) small
w.r.t.\ the left member. This would yield the same result without any
restriction on $\mu_k$, at the expense of a somewhat longer proof.

One problem in polling models with $1$-limited service strategy is
that the rank of the systems of equations that we can write is
$\lfloor n/2\rfloor$, which means that $\dpar{z_m}F(\e)$ cannot be
computed. Fortunately, under the following assumption, we are able to
compute the most important value, that is $\dpar{z_m}{F_m}(\e)$.

\begin{ass}\label{ass:miroir}
The routing matrices are such that
\[
 [I-P][I-\tr{\wt P}]=[I-\tr P][I-\wt P].
\]
\end{ass}

It is not easy to describe all models satisfying~\ref{ass:miroir}.
The simplest one is the classical Markov polling obtained when $P=\wt
P$, for which the results of Theorem~\ref{thm:main} thereafter are
greatly simplified. Another admissible model is when the routing
probabilities depend only on the absolute distance between stations,
that is when $P=\tr P$ and $\wt P=\tr{\wt P}$.

\begin{thm}\label{thm:main}
Assume that~\ref{ass:sym},~\ref{ass:irred} and~\ref{ass:miroir} hold.
Then, for any station $m$, the stationary probability that the server
finds station $m$ empty is given by
\begin{equation}\label{eq:Fte}
P(X_m=0\mid S=m)={1-N\alpha\over1-N\alpha+N\t\alpha}\;.
\end{equation}

Moreover the mean number of customers found by the server when it
arrives at station $m$ is obtained by means of the following formula
\begin{eqnarray}\label{eq:dFe}
\EE[X_m\mid S=m]
&=& N\bar\alpha+{N\t\alpha\bar\alpha\over 1-N\alpha}
        \sum_{l=1}^{N-1}{1\over 1-\t\mu_l}\nonumber\\
& &\mbox{}+{1-N\alpha+N\t\alpha\over 1-N\alpha}{N\over2}\bar\alpha^\2_N
   +{N\alpha-N\t\alpha\over 1-N\alpha}
      {1\over2}\sum_{d=1}^N \bar\alpha^\2_d\nonumber\\
& &\mbox{}-{N\t\alpha\over 1-N\alpha}
      \sum_{l=1}^{N-1}{\mu_l-\t\mu_l\over1-\t\mu_l}
         {1\over2}\sum_{d=1}^N\bar\alpha^\2_d\omega_l^d\;.
\end{eqnarray}

Finally, the mean number of customers present at arbitrary station $m$
at a polling instant can be expressed as
\begin{eqnarray}\label{eq:sumdFe}
\EE X_m
&=& \bar\alpha
    +{\t\alpha\over 1-N\alpha}
        \sum_{l=1}^{N-1}{1\over 1-\t\mu_l}\nonumber\\
& &\mbox{}+{1-N\alpha+N\t\alpha\over 1-N\alpha}{N\over2}\bar\alpha^\2_N
   +{N\alpha-N\t\alpha\over 1-N\alpha}
      {1\over2}\sum_{d=1}^N \bar\alpha^\2_d\nonumber\\
& &\mbox{}-{1-N\alpha+N\t\alpha\over 1-N\alpha}
      \sum_{l=1}^{N-1}{\mu_l-\t\mu_l\over1-\t\mu_l}
         {1\over2}\sum_{d=1}^N\bar\alpha^\2_d\omega_l^d\;.
\end{eqnarray}
\end{thm}

\begin{proof}{} 
See Appendix~\ref{app:proof}. 
\end{proof}

\section{The mean waiting time of a customer}\label{sec:wait}

Theorem~\ref{thm:main} is the key to compute of the waiting time of an
arbitrary customer. In this section, we focus on the continuous-time
case with compound Poisson arrivals and further assume that the
arrivals at each station are independent. While Theorem~\ref{thm:main}
is valid in discrete time or with correlated arrivals, we do not give
any formula in these cases, since the notation would become too
involved.

The notation is as follows: customers arrive in i.i.d.\ batches of
mean length $b$ and second moment $b^\2$ at the instants of a Poisson
stream with intensity $\hat\lambda$ at each station. The intensity of
the arrival process at each queue is $\lambda=\hat\lambda b$. The time
between two polling instants when the first queue is not empty (resp.\
empty) has mean and second moment $\tau$ and $\tau^\2$ (resp.\
$\t\tau$ and $\t\tau^\2$). We must keep in mind that in general $\tau$
depends on $P$ like $\alpha_{i;q}$ in Equation~(\ref{eq:defaliq}): if
$\tau_d$ is the mean time between two consecutive polling instants at
stations $i$ and $i+d$, then we have
$\tau=p_1\tau_1+\cdots+p_N\tau_N$. However, in most symmetrical
polling models that have been previously analyzed, the switchover
times are equal and $\tau$ does not depend on $P$. We will say in this
case that the stations are {\em equidistant}, which is obviously not
the case in our original taxicab problem. The same holds for
$\tau^\2$, $\t\tau$ and $\t\tau^\2$.

\begin{thm}\label{thm:EW}
Assume that~\ref{ass:sym},~\ref{ass:irred} and~\ref{ass:miroir} hold.
Then the mean stationary waiting time of a customer is given by the
relation
\begin{eqnarray}\label{eq:EW}
\EE[W]
&=& {\t\tau\over 1-N\lambda\tau}
      \sum_{l=1}^{N-1}\biggl[{1\over 1-\t\mu_l}\biggr]
    +{N\lambda \tau^\2\over 2(1-N\lambda \tau)}
      +{\t\tau^\2\over2\t\tau}\nonumber\\
\bigskip\nonumber\\
& &\mbox{}+ {b^\2-b\over 2b}\Biggl[
      {\tau+(N-1)\t\tau\over 1-N\lambda\tau}
     -{\t\tau\over 1-N\lambda\tau}
\sum_{l=1}^{N-1}{\mu_l-\t\mu_l\over 1-\t\mu_l}\Biggr].
\end{eqnarray}
\end{thm}

\begin{proof}{}
See Appendix~\ref{app:proof}.
\end{proof}

One surprising consequence of this formula is that when the
arrival process is Poisson and the stations are {\em equidistant}, the
mean waiting time does not depend on the routing $P$. In other words,
it does not depend on where the server goes after serving a customer.

Beside giving the mean waiting time for a customer in our taxicab
model, Equation~(\ref{eq:EW}) contains in fact many known formulas for
waiting times corresponding to various service and polling
strategies. In the next subsections, we give some of these
applications.

A notation closer to the classical queuing theory notation is
necessary to make the link with known results. Let $w$ (resp.\ $\t w$)
be the mean walking time---or switchover time---from a non-empty
(resp.\ empty) station and let $\sigma$ be the mean service time
required by the customers. We denote by $w^\2$, $\t w^\2$ and $\sigma^\2$
the associated second moments.

\subsection{A state-independent Markovian polling model}

The first possible application of our model is the classical state
independent polling model with $1$-limited service strategy. In this
model, we have $P=\wt P$, $\tau=w+\sigma$, $\t\tau=w$
and~(\ref{eq:EW}) becomes
\begin{eqnarray*}
\EE[W]
&=& {w\over1-N\lambda(w+\sigma)}
      \sum_{l=1}^{N-1}{1\over1-\mu_l}
    +{N\lambda(w^\2+2w\sigma+\sigma^\2)
         \over 2(1-N\lambda(w+\sigma))}
    +{w^\2\over 2w}\\
& &\mbox{}-{b^\2-b\over 2b}{Nw+\sigma\over1-N\lambda(w+\sigma)}.
\end{eqnarray*}

This formula is valid for all symmetric polling models with Markovian
polling. It is interesting to note that, in the case of equidistant
stations, there is only one term of the formula which depends on the
routing. To compute this term, we remark that, if ${\cal P}(x)$ is the
characteristic polynomial of $P$,
\[
\sum_{l=1}^{N-1}{1\over 1-\mu_l}
\;=\; \sum_{l=1}^{N-1}{d\over dx}\log(x-\mu_l)\bigg|_{x=1}
\;=\; {d\over dx}\log\Bigl[{{\cal P}(x)\over x-1}\Bigr]\bigg|_{x=1}.
\]

The case of the cyclic polling is obtained by taking $p_1=1$ and
$p_d=0$ for $d\neq1$. Then the eigenvalues of $P$ are
$\omega_1,\ldots,\omega_N$ and ${\cal P}(x)=x^N-1$. So in this case,
we have
\[ \sum_{l=1}^{N-1}{1\over 1-\mu_l}={N-1\over2}. \]

Another classical polling strategy is the random polling with
$p_d=1/N$ for all $d$. In this case, we have $\mu_l=0$ for $l<N$ and
\[ \sum_{l=1}^{N-1}{1\over 1-\mu_l}=N-1. \]

The comparison between formulas for cyclic polling and random polling
shows that the mean waiting time is smaller in the cyclic case. In
fact this property can be generalized to all Markovian polling
models.

\begin{lem}\label{lem:min}
Among all Markovian polling strategies for symmetric $1$-limited
polling models with {\em equidistant\/} queues, the cyclic polling
strategy minimizes the waiting time of the customers.
\end{lem}
\begin{proof}{}
This result is very easy to prove from Equation~(\ref{eq:EW}). We have
\begin{eqnarray*}
\sum_{l=1}^{N-1}{1\over 1-\mu_l}
&=& \sum_{l=1}^{N-1}\Re\Bigl({1\over 1-\mu_l}\Bigr)\\
&=& \sum_{l=1}^{N-1}{1-\Re(\mu_l)\over 2(1-\Re(\mu_l))+|\mu_l|^2-1},
\end{eqnarray*}
where $\Re(z)$ is the real part of $z$. Since $|\mu_l|\leq1$ and
$\Re(\mu_l)\leq1$, we find that
\[  \sum_{l=1}^{N-1}{1\over 1-\mu_l}\geq {N-1\over 2}. \]

The bound is attained if and only if for all $l$, $|\mu_l|=1$.  Since
$\mu_l$ is the center of gravity of $\omega_l^1,\ldots,\omega_l^N$
with weights $p_1,\ldots,p_N$, this is only possible when for some
$s\in\S$ we have $p_s=1$ and $p_d=0$ for $d\neq s$. When $s$ is not a
divider of $N$ (to ensure that the routing matrix is irreducible),
this is equivalent to cyclic polling.
\end{proof}

This property also holds for discrete time systems and systems with
correlated arrivals. The computation would be more difficult in the
case where the stations are not equidistant.

\subsection{The exhaustive and Bernoulli service strategies}

An interesting application of our model is to show that some service
strategies can be obtained by a proper choice of polling strategy. In
this subsection, we show how we can apply our model to exhaustive and
Bernoulli service strategies. 

Consider the case where the stations are {\em equidistant}, except
that the switchover time from a station to itself after a service is
zero.  Then, if we denote $p_N=1-\pi$, we have $\tau=\pi w+\sigma$,
$\tau^\2=\pi w^\2+2\pi w\sigma+\sigma^\2$, $\t\tau=w$ and
$\t\tau^\2=w^\2$. If we assume for the sake of simplicity that the
arrival process is Poisson,~(\ref{eq:EW}) yields
\[
\EE[W] = {w\over1-N\lambda(\pi w+\sigma)}
      \sum_{l=1}^{N-1}{1\over1-\t\mu_l}
    +{N\lambda(\pi w^\2+2\pi w\sigma+\sigma^\2)
         \over 2(1-N\lambda(\pi w+\sigma))}
    +{w^\2\over 2w}.
\]

The value of $\EE[W]$ depends on $P$ only through the value of $p_N$.
One known model that is described by this equation is the Bernoulli
strategy of Servi~\cite{Serv:1}: when the server has served a
customer, it quits the queue with probability $\pi$ and continues to
serve it with probability $1-\pi$. In the original model, the server
polls the queues in cyclic order. Here, we have a Markovian Bernoulli
polling model if we take $P=(1-\pi)I+\pi\wt P$. It is easy to check
that this choice of $P$ satisfies~\ref{ass:miroir}. Moreover, $\EE[W]$
is an increasing function of $\pi$ that is minimal when $\pi=0$. This
case corresponds to an exhaustive service strategy with Markovian
routing $\wt P$: the server polls the same queue until it is empty and
then moves to another queue using the routing matrix $\wt P$. In this
case, Equation~(\ref{eq:EW}) simply becomes
\[
 \EE[W] \;=\; {w\over1-N\lambda\sigma}
      \sum_{l=1}^{N-1}{1\over1-\t\mu_l}
    +{N\lambda\sigma^\2
         \over 2(1-N\lambda\sigma)}
    +{w^\2\over 2w}.
\]

Using Lemma~\ref{lem:min}, we find that the above expression is 
minimal when $\wt P$ describes a cyclic polling scheme. Further 
comparisons between different polling strategies can be found 
in Levy, Sidi and Boxma~\cite{LevSidBox:1}.

\appendix
\section*{Appendix}
\def\thesection{A}
\subsection{Second vector field}\label{app:induc}

This appendix contains some basic definitions and results from the
theory of dynamical systems. These definitions have been adapted for
the Markov chain $\L$ defined in Section~\ref{sec:model}. We refer the
reader to~\cite{MalMen:1} and~\cite{FayMalMen:1} for a more complete
treatment of the subject.

\medskip\noindent
{\bf Faces}. For any $\l\subset\{1,\ldots,N\}$, define
$B^\l\subset\RR^N_+$ as
\[
  B^\l\egaldef \{(x_1,\ldots,x_n): x_i>0 \Leftrightarrow i\in\l\}.
\]

$B^\l$ is the {\em face\/} of $\RR^N_+$ associated to $\l$. Whenever
not ambiguous, the face $B^\l$ will be identified with $\l$. 

\medskip\noindent
{\bf Induced chains}. For any $\l \neq \S$ we choose an arbitrary
point $\vecteur a\in B^\l\cap \ZZ^N_+$ and draw a plane $C^\l$ of
dimension $N-|\l|$, perpendicular to $B^\l$ and containing $\vecteur
a$. We define the {\em induced Markov chain\/} $\L^\l$, with state
space $\S\times (C^\l\cap\ZZ^N_+)$ (by an obvious abuse in the
notation, we shall write most of the time $\S\times C^\l$) and
transition probabilities
\[
{}_\l p_{(s,\x)(t,\y)}=
  p_{(s,\x)(t,\y)}+\sum_{\z \neq \y}p_{(s,\x)(t,\z)}, 
  \ \forall \x,\y \in C^{\l},\ s,t\in\S,
\]
where the summation is performed over all $\z\in\ZZ^N_+$, such
that the straight line connecting $\z$ and $\y$ is
perpendicular to $C^\l$. It is important to note that this
construction does not depend on $\vecteur a$.

\begin{ass}
For any $\l$ the chain $\L^\l$ is irreducible and aperiodic.  $\l$ is
called {\em ergodic\/} (resp.\ {\em non-ergodic, transient\/})
according as ${\cal L}^\l$ is ergodic (resp.\ non-ergodic, transient).
\end{ass}

For an ergodic $\L^\l$, let $\pi^\l(s,\x),\;(s,\x)\in \S\times C^\l$
be its stationary transition probabilities. Moreover, let $\pi^\l(s)$
(resp.\ $\pi^\l(s,x_s>0)$) be the stationary probability that the
server is at station $s$ (resp.\ and finds it nonempty).

\medskip\noindent
{\bf First vector field}. For any $(s,\x)\in\S\times \ZZ^N_+$, the
{\em first vector field\/} is simply the mean drift of the random walk
$\L$ at point $(s,\x)$:
\[
 \vecteur M(s,\x)\egaldef 
    \!\!\!\sum_{r\in\S,\;\y\in\ZZ_+^N}\!\!\!(\y-\x)
       P\Bigl[X(n+1)=\y,S(n+1)=r \,\Big|\, X(n)=\x, S(n)=s\Bigr].
\]

Let $\v^\l\egaldef(v^\l_1,\ldots,v^\l_N)$ such that
\begin{eqnarray*}
v^\l_i & = & 0,\; i\not \in \l,\\
v^\l_i & = & \sum_{(s,\x)\in \S\times C^\l} \; \pi^\l(s,\x)\vecteur M_i(s,\x),
                      \; i\in \l.
\end{eqnarray*}

Intuitively, one can imagine that the random walk starts from a point
which is close to $\l$, but sufficiently far from all other faces
$B^{\l'}$, with $\l\not\subset\l'$. After some time (sufficiently
long, but less than the minimal distance from the above mentioned
$B^{\l'}$), the stationary regime in the induced chain will be
installed. In this regime, one can ask about the mean drift along
$\l$: it is defined exactly by $\v^\l$. For $\l =\S$, we call $\l$ {\it
ergodic}, by definition, and put
\[
  \v^{\{1,\ldots,N\}}\equiv 
    \sum_{s\in\S}\pi^{\{1,\ldots,N\}}(s)\vecteur M(s,\x), 
  \; \x\in  B^{\{1,\ldots,N\}}.
\] 

From now on, when speaking about the {\em components\/} of $\v^\l$, we
mean the components $v^\l_i$ with $i\in\l$.

\begin{ass}
$v^\l_i\neq 0$, for each $i\in\l$.
\end{ass}

\noindent
{\bf Ingoing, outgoing and neutral faces.} Let us fix $\l ,\l_1$, so
that $\l\supset\l_1, \l \neq \l_1$, that is to say
$\overline{B}^\l\supset B^{\l _1}$ ($\overline{B}^\l$ is the closure
of $B^\l$). Let $B^\l$ be ergodic.  Thus $\v^\l$ is well defined. There
are three possibilities for the direction of $\v^\l$ w.r.t. $B^{\l
_1}$.  We say that $B^\l$ is an {\em ingoing\/} (resp.\ {\em
outgoing\/}) {\em face\/} for $B^{\l_1}$, if all the coordinates
$v^\l_i$ for $i\in \l\setminus\l_1$ are negative (resp.\
positive). Otherwise we say that $B^\l$ is {\em neutral}.  As an
example we give simple sufficient criteria for a face to be ergodic.

\medskip\noindent
{\bf The second vector field}. To any point $\x\in \RR^N_+$, we
assign a vector $v(\x)$ and call this function the {\em second vector field}.
It can be multivalued on some non-ergodic faces. We put, for ergodic
faces $B^\l$,
\[
  \v(\x)\equiv \v^\l, \; \x\in B^\l.
\]

If $B^{\l _1}$ is non-ergodic, then at any point $\x\in
B^{\l _1}$, $\v(\x)$ takes all values $\v^\l$ for
which $B^\l$ is an outgoing face with respect to
$B^{\l _1}$. In other words, for $\x$ belonging to non-ergodic faces,
with $\|\x\|$ sufficiently large,
\[
  \x+\v(\x)\;\in \RR^N_+\, ,
\]
for any value $\v(\x)$. If there is no such vector, we put $\v(\x)=0$,
for $\x\in B^{\l _1}$. Points $\x\in {\bf R}^N_+$, where $\v(\x)$ is
more than one-valued, are called branch points.

There are few interesting examples, for which only the first vector
field suffices to obtain ergodicity conditions for the random walk of
interest, but it is nevertheless the case for Jackson networks. In
general, the second vector field must be introduced as shown in the
following theorem, taken from~\cite{MalMen:1} (and extended
in~\cite{FayZam:1} to the case of upward unbounded jumps).

\begin{thm}\label{thm:FayZam}
Assume that there exists a nonnegative function
$f\,:\,\RR^N_+\to\RR_+$ such that
\begin{enumerate}
\item for some constant $c>0$
\[
 f(\x)-f(\y)\leq c\|\x-\y\|;
\]
\item there exists $\delta>0$, $p>0$ such that for all $x\in B^\l$
\[
 f(\x+v(\x))-f(\x)\leq -\delta.
\]
\end{enumerate}

Then the Markov chain $\L$ is ergodic.
\end{thm}

\subsection{Some simple results about circulant matrices}\label{app:circu}

In this appendix, we recall some well-known properties of the
circulant matrices which are used throughout this paper. These
properties are given without proof, since they can easily be verified
at hand.  Throughout the paper, we take the convention that every
subscript less than 1 or greater than $N$ should be shifted into the
correct range.

Let $(\e_1,\ldots,\e_N)$ denote the canonical basis of $\CC^N$ and
let $\e=(1,\ldots,1)$. Moreover, we define $\omega_k\egaldef
\exp(2\hat\iota\pi k/N)$, with $\hat\iota^2=-1$.

\begin{defi}
A circulant matrix $M$ of size $N$ is a matrix whose coefficients
$m_{i,j}$ verify the relation:
\[
 m_{i+k,j+k}=m_{i,j}, \mbox{\rm \ for all } 1\leq i,j,k\leq N.
\]
\end{defi}

The following Lemma summarizes some key properties of these matrices:

\begin{lem}\label{lem:circu}
Let $M$ be a circulant matrix of the form:
\[
M=\left(
\begin{array}{cccc}
 m_N    & m_{N-1} & \cdots  & m_1   \\
 m_1    & m_N     & \cdots  & m_2   \\
 \vdots &         & \ddots & \vdots\\
 m_{N-1}& m_{N-2} & \cdots  & m_N   
\end{array}
\right)
\]

Then for each $1\leq k\leq N$, the vector $\v_k\egaldef\sum_{i=1}^N
\omega_k^{-i}\e_i$ is an eigenvector of the matrix $M$ with eigenvalue
$\sum_{i=1}^N\omega_k^im_i$. Moreover, $(\v_1,\ldots,\v_N)$ is an
orthogonal basis of $\CC^N$ in which $\e_1,\ldots,\e_N$ can be
expressed as
\[
\e_i={1\over N}\sum_{k=1}^N \omega_k^i\v_k.
\]
\end{lem}

One important feature of circulant matrices is that they share the
same basis of eigenvectors. Note that, with the notation given before
Lemma~\ref{lem:circu}, we have $\v_N=\e$. Circulant matrices enjoy
other properties that are not used here: for example, the product
of two circulant matrices is a circulant matrix and this product
commutes.

\subsection{Proofs of the results of Sections \protect\ref{sec:sym} 
         and \protect\ref{sec:wait}}
\label{app:proof}

The derivations in this section are essentially analytic. 

\begin{proof}{ of Lemma~\ref{lem:ut}}
This proof uses specific properties of circulant matrices to perform a
fine analysis of the behavior of $A(\z)$, $\wt A(\z)$ and
$\Delta(\z)$ in the neighborhood of $\z=\e$. As pointed out later in
the proof, we study these matrices only for $\z\in{\cal D}^N$, thus
avoiding any analytical continuation. The basic relation used
thereafter is a simple consequence of the definitions of $A(\z)$ and
$\Delta(\z)$:
\begin{equation}\label{eq:Avk}
 A\Delta(\z)\v_k=\mu_k(\z)\v_k
          +\sum_{q=1}^N({1\over z_q}-1)\omega_k^{-q}A_{\cdot q}(\z),
\end{equation}
where $A_{\cdot q}(\z)$ stands for the $q$-th column of $A(\z)$ and
can be written as
\begin{eqnarray*}
 A_{\cdot q}(\z)
&=& \sum_{i=1}^N p_{i-q} a_{i-q}(\z)\e_i\\
&=& \sum_{i=1}^N p_{i-q} a_{i-q}(\z){1\over N}\sum_{l=1}^N \omega_l^i\v_l\\
&=& \sum_{l=1}^N \omega_l^q\mu_l(\z)\v_l.
\end{eqnarray*}

So, if we define 
\[
\epsilon_l(\z)\egaldef {1\over N}\sum_{q=1}^N({1\over z_q}-1)\omega_l^q,
\]
(\ref{eq:Avk}) can be rewritten as
\begin{eqnarray}
[I-A\Delta(\z)]\v_k
 &=& (1-\mu_k(\z))\v_k
     -\sum_{l=1}^N\epsilon_{l-k}(\z)\mu_l(\z)\,\v_l\label{eq:ImAvk}\\
\mbox{}[\wt A(\z)-A\Delta(\z)]\v_k
 &=& (\t\mu_k(\z)-\mu_k(\z))\v_k
     -\sum_{l=1}^N\epsilon_{l-k}(\z)\mu_l(\z)\,\v_l\label{eq:AtmAvk}
\end{eqnarray}

Since $\tr\e[I-A\Delta(\e)]=\tr\e[\wt A(\e)-A\Delta(\e)]=\vecteur0$,
we define, for any set of arbitrary complex numbers
$(c_1,\ldots,c_{N-1})$, the vector $\u_t$ as follows:
\[
  \u_t\egaldef {1\over N}\e+{t\over N}\sum_{l=1}^{N-1}c_l\v_l.
\]

With this definition, we have 
\begin{eqnarray}
\tr{\u_t}[I-A\Delta(\z)]\e
&=& {\tr\e\over N}[I-A\Delta(\z)]\e+
    {t\over N}\sum_{l=1}^{N-1}c_l\tr{\v_l}[I-A\Delta(\z)]\e\nonumber\\
&=& 1-\mu_N(\z)-\epsilon_0(\z)\mu_N(\z)
    -t\sum_{l=1}^{N-1}c_l\epsilon_l(\z)\mu_l(\z)\label{eq:tue}.
\end{eqnarray}

When $t$ is small, we have the relations:
\begin{eqnarray}
\mu_N(\z_t)
 &=& 1+t\alpha\sum_{q=1}^N \dot z_q
     +{t^2\over2}\sum_{q,r=1}^N\alpha^\2_{q,r}\dot z_q\dot z_r\nonumber\\
 & &\mbox{}+{t^2\alpha\over 2}\sum_{q=1}^N\ddot z_q+o(t^2)\label{eq:dlmu}\\
\t\mu_N(\z_t)
 &=& 1+t\t\alpha\sum_{q=1}^N \dot z_q
     +{t^2\over2}\sum_{q,r=1}^N\t\alpha^\2_{q,r}\dot z_q\dot z_r\nonumber\\
 & & \mbox{}+{t^2\t\alpha\over 2}\sum_{q=1}^N\ddot z_q+o(t^2)
                                                            \label{eq:dltmu}\\
\epsilon_l(\z_t)
 &=& -{t\over N}\sum_{q=1}^N \dot z_q\omega_l^q
     +{t^2\over N}\sum_{q=1}^N\dot z_q^2\omega_l^q\nonumber\\
 & & \mbox{}-{t^2\over2N}\sum_{q=1}^N\ddot z_q\omega_l^q+o(t^2).
                                                            \label{eq:dleps}
\end{eqnarray}

The formulas given in Appendix~\ref{app:circu} allow to relate $\dot
z_1,\ldots,\dot z_N$ to $\dot\zeta_1,\ldots,\dot\zeta_N$:
\[
 \dot\zeta_k={1\over N}\sum_{q=1}^N\dot z_q\omega_k^q,\qquad
 \dot z_q=\sum_{l=1}^N\dot\zeta_q\omega_l^{-q},\qquad
 \sum_{q=1}^N\dot z_q\dot z_{q+d}
   =N\sum_{l=1}^N\dot\zeta_l\dot\zeta_{N-l}\omega_l^d.
\]

Using the remark after the definition of $\alpha^\2_{q,r}$ and
$\t\alpha^\2_{q,r}$, we find that:
\[
 \sum_{q,r=1}^N\alpha^\2_{q,r}\dot z_q\dot z_r
 =\sum_{d=1}^N\alpha^\2_d\sum_{q=1}^N\dot z_q\dot z_{q+d}
 =N\sum_{l=1}^N\dot\zeta_l\dot\zeta_{N-l}\sum_{d=1}^N\alpha^\2_d\omega_l^d.
\]

Applying~(\ref{eq:dlmu}) and~(\ref{eq:dleps}) to~(\ref{eq:tue}), we
see that the first-order term in the expression~(\ref{eq:tue}) is
$(\alpha-1/N)t[\dot z_1+\cdots+\dot z_N]$. So a necessary condition to
have a formula like~(\ref{eq:dltue}) is $\dot z_1+\cdots+\dot
z_N=0$. Using this relation, we have
\begin{eqnarray}
\tr{\u_t}[I-A\Delta(\z_t)]\e
 &=&  t^2\sum_{l=1}^{N-1}c_l\dot\zeta_l\mu_l
  -t^2\sum_{l=1}^N\dot\zeta_l\dot\zeta_{N-l}\label{eq:dltue1}\nonumber\\
 & & \mbox{}-t^2{N\over2}\sum_{l=1}^N\dot\zeta_l\dot\zeta_{N-l}
                            \sum_{d=1}^N\alpha^\2_d\omega_l^d\nonumber\\
 & &\mbox{}+t^2{1-N\alpha\over2N}\sum_{q=1}^N\ddot z_q+o(t^2).
\end{eqnarray}

We have to check that it is possible to have $\z_t\in{\cal D}^N$ and
$\dot z_1+\cdots+\dot z_N=0$. One easy way to satisfy these
constraints is to ensure that all coordinates of $\z_t$ arrive to $1$
tangentially to the unit circle as $t$ goes to $0$. This is the case
when, for any $q$, $\dot z_q$ is an imaginary number and $\ddot
z_q<0$. Moreover, for $k<N$,
\begin{eqnarray}\label{eq:tuvk}
\tr{\u_t}[I-A\Delta(\z_t)]\v_k
&=& -\epsilon_{N-k}(\z_t)\mu_N(\z_t)+tc_k(1-\mu_k(\z_t))\nonumber\\
& & \mbox{}-t\sum_{l=1}^{N-1}c_l\epsilon_{l-k}(\z_t)\mu_l(\z_t)\nonumber\\
&=& t\dot\zeta_{N-k}+tc_k(1-\mu_k)+o(t).
\end{eqnarray}

This shows that we get equations~(\ref{eq:dltuvk})
and~(\ref{eq:dltue}) from~(\ref{eq:dltue1}) and~(\ref{eq:tuvk}) if we
take
\[
  c_k=-{\dot\zeta_{N-k}\over1-\mu_k}.
\]

With this choice of $c_1,\ldots,c_{N-1}$, we obtain
equations~(\ref{eq:dltuAte}) and~(\ref{eq:dltuAtvk}) in exactly the
same way.
\end{proof}

\begin{proof}{ of Theorem~\ref{thm:main}}
The idea of this proof is to apply the results of Lemma~\ref{lem:ut}
to the equation
\begin{equation}\label{eq:tumatrix}
\tr{\u_t}[I-A\Delta(\z_t)]F(\z_t)
  =\tr{\u_t}[\wt A(\z_t)-A\Delta(\z_t)]\wt F(\z_t).
\end{equation}

Define $\z_t$ as in Lemma~\ref{lem:ut}. Then

\[
 F(\z_t)=F(\e)+t\sum_{r=1}^N\dot z_r\d{z_r}F(\e)+\o(t).
\]

Moreover, Equation~(\ref{eq:symF}) implies that
$\dpar{z_r}{F_i}(\e)=\dpar{z_m}{F_{i-r+m}}(\e)$ for any fixed
$m\in\S$ and
\begin{eqnarray*}
 \d{z_r}F(\e)
&=& \sum_{i=1}^N\d{z_r}{F_i}(\e)\e_i\\
&=& \sum_{i=1}^N\Bigl[\d{z_m}{F_{i-r+m}}(\e)
      {1\over N}\sum_{k=1}^N\omega_k^i\v_k\Bigr]\\
&=& \sum_{k=1}^N\omega_k^{r-m}
               \Bigl[{1\over N}\sum_{i=1}^N\d{z_m}{F_{i-r+m}}(\e)
                                       \omega_k^{i-r+m}\Bigr]\v_k\\
&=& \sum_{k=1}^N\omega_k^{r-m}\d{z_m}{\phi_k}(\e)\v_k,
\end{eqnarray*}
where $(\phi_1(\z),\ldots,\phi_N(\z))$, defined as
\[
 \phi_k(\z)\egaldef{1\over N}\sum_{i=1}^NF_i(\z)\omega_k^i,
\]
are the coordinates of $F(\z)$ in the basis
$(\v_1,\ldots,\v_N)$. Finally,
\begin{equation}\label{eq:Fzt}
 F(\z_t)=F(\e)+tN\sum_{k=1}^N\dot\zeta_k
    \omega_k^{-m}\d{z_m}{\phi_k}(\e)\v_k+\o(t).
\end{equation}

With a similar definition for $\t\phi(\z)$,
\begin{equation}\label{eq:tFzt}
 \wt F(\z_t)=\wt F(\e)+tN\sum_{k=1}^N\dot\zeta_k
    \omega_k^{-m}\d{z_m}{\t\phi_k}(\e)\v_k+\o(t).
\end{equation}

We now apply Lemma~\ref{lem:ut} and Equations~(\ref{eq:Fzt})
and~(\ref{eq:tFzt}) to Equation~(\ref{eq:tumatrix}) and use the fact
that, by symmetry, $F(\e)=F_1(\e)\e$ and $\wt F(\e)=\wt F_1(\e)\e$

\begin{eqnarray*}
\lefteqn{t^2F_1(\e)\Biggl\{-\sum_{l=1}^{N-1}
  \dot\zeta_l\dot\zeta_{N-l}\biggl[
   {1\over 1-\mu_l}+{N\over2}\sum_{d=1}^N\alpha^\2_d\omega_l^d\biggr]
  +{1-N\alpha\over2N}\sum_{q=1}^N\ddot z_q\Biggr\}}\qqqq\\
 &=& t^2\wt F_1(\e)\Biggl\{-\sum_{l=1}^{N-1}
  \dot\zeta_l\dot\zeta_{N-l}\biggl[
  {1\over 1-\mu_l}
  +{N\over2}\sum_{d=1}^N(\alpha^\2_d-\t\alpha^\2_d)\omega_k^d\biggr]\\
 & &\qqqq+{1-N\alpha+N\t\alpha\over2N}\sum_{q=1}^N\ddot z_q\Biggr\}\\
 & &\mbox{}+t^2N\sum_{k=1}^{N-1}\dot\zeta_k\dot\zeta_{N-k}\omega_k^{-m}
     {1-\t\mu_k\over 1-\mu_k}\d{z_m}{\wt\phi_k}(\e)+o(t^2).
\end{eqnarray*}

As $\ddz$ can be chosen freely, we can derive a first equality from
this equation, namely

\[
 F_1(\e){1-N\alpha\over2N}=\wt F_1(\e){1-N\alpha+N\t\alpha\over2N}.
\]

Since $P(X_m=0\mid S=m)=\wt F_1(\e)/F_1(\e)$, this yields
Equation~(\ref{eq:Fte}).  Taking in account~(\ref{eq:bal2}), we find

\begin{eqnarray}\label{eq:totalphit}
\lefteqn{ \sum_{k=1}^{N-1}\dot\zeta_k\dot\zeta_{N-k}
 \omega_k^{-m}{1-\t\mu_k\over 1-\mu_k}N\d{z_m}{\wt\phi_k}(\e)}\qqqq\nonumber\\
&=& -\sum_{l=1}^{N-1}\dot\zeta_l\dot\zeta_{N-l}\biggl[
  {F_1(\e)-\wt F_1(\e)\over 1-\mu_l}
  +{1\over2}\sum_{d=1}^N\bar\alpha^\2_d\omega_l^d\biggr].
\end{eqnarray}

This equation contains in fact a system of linear equations that can
be built by choosing $\dz$. The problem is that in general the rank of
this system is only $\lfloor n/2\rfloor$. However, since $\mu_k$
(resp.\ $\mu_{N-k}$, $\t\mu_k$, $\t\mu_{N-k}$) is the eigenvalue of
$\tr P$ (resp.\ $P$, $\tr{\wt P}$, $\wt P$) associated to the
eigenvector $\v_k$,~\ref{ass:miroir} implies that, for
$1\leq k<N$,
\[
 {1-\mu_k\over1-\t\mu_k}={1-\mu_{N-k}\over1-\t\mu_{N-k}}\in \RR.
\]

As noted in the proof of Lemma~\ref{lem:ut}, $\dot z_1,\ldots,\dot
z_N$ are imaginary numbers and $\dot\zeta_k\dot\zeta_{N-k}$ is a
negative real number for any $k<N$. Hence we can choose
in~(\ref{eq:totalphit})
\[
 \dot\zeta_k\dot\zeta_{N-k}=
  -\min\Bigl[{1-\mu_k\over 1-\t\mu_k},0\Bigr].
\]

The combination of the resulting equation with a similar equation
 containing only the terms where $(1-\mu_k)/(1-\t\mu_k)<0$ yields
\begin{eqnarray}\label{eq:totalphit2}
\lefteqn{ \sum_{k=1}^{N-1}\omega_k^{-m}N\d{z_m}{\wt\phi_k}(\e)}\qqqq\nonumber\\
&=& -\sum_{l=1}^{N-1}
  {1-\mu_k\over 1-\t\mu_k}\biggl[
  {F_1(\e)-\wt F_1(\e)\over 1-\mu_l}
  +{1\over2}\sum_{d=1}^N\bar\alpha^\2_d\omega_l^d\biggr].
\end{eqnarray}

We know that, by definition, $\wt{F}_m(\z)$ does not depend on $z_m$.
Hence,
\[
 \d{z_m}{\wt F_m}(\e)=\sum_{k=1}^N\omega_k^{-m}\d{z_m}{\wt\phi_k}(\e)=0,
\]
and, using Equation~(\ref{eq:totalphit2})
\begin{equation}\label{eq:dphiN}
N\d{z_m}{\t\phi_N}(\e)
 = \sum_{l=1}^{N-1}\Biggl[
  {F_1(\e)-\wt F_1(\e)\over1-\t\mu_l}
  +{1-\mu_l\over1-\t\mu_l}{1\over2}\sum_{d=1}^N
          \bar\alpha^\2_d\omega_l^d\Biggr].
\end{equation}

The last step of the demonstration is to get $\dpar{z_m}{F_m}(\e)$ from
these results. This can be done as in Lemma~\ref{lem:ut}, but with
$\u_t=\e/N$ and $\z_t$ chosen differently. Using
Equations~(\ref{eq:ImAvk}),~(\ref{eq:AtmAvk})
and~(\ref{eq:dlmu})--(\ref{eq:dleps}) with the same notation as in
Lemma~\ref{lem:ut}, we have for $k<N$
\begin{eqnarray*}
{\tr\e\over N}[I-A\Delta(\z_t)]\e 
  &=& 1-\mu_N(\z_t)-\epsilon_0(\z_t)\mu_N(\z_t)\\
  &=& t(1-N\alpha)\dot\zeta_N
      -t^2\sum_{l=1}^N\dot\zeta_l\dot\zeta_{N-l}+t^2N\alpha\dot\zeta_N^2\\ 
  & &\mbox{}-t^2{N\over2}\sum_{l=1}^N\dot\zeta_l\dot\zeta_{N-l}
           \sum_{d=1}^N\alpha^\2_d\omega_l^d+o(t^2)\\
{\tr\e\over N}[I-A\Delta(\z_t)]\v_k
  &=& -\epsilon_{N-k}(\z_t)\mu_N(\z_t) \;=\; t\dot\zeta_{N-k}+o(t)\\
{\tr\e\over N}[\wt A(\z_t)-A\Delta(\z_t)]\e 
  &=& \t\mu_N(\z_t)-\mu_N(\z_t)-\epsilon_0(\z_t)\mu_N(\z_t)\\
  &=& t(1-N\alpha+N\t\alpha)\dot\zeta_N
      -t^2\sum_{l=1}^N\dot\zeta_l\dot\zeta_{N-l}+t^2N\alpha\dot\zeta_N^2\\
  & &\mbox{} +t^2{N\over2}\sum_{l=1}^N\dot\zeta_l\dot\zeta_{N-l}
           \sum_{d=1}^N[\t\alpha^\2_d-\alpha^\2_d]\omega_l^d+o(t^2)\\
{\tr\e\over N}[\wt A(\z_t)-A\Delta(\z_t)]\v_k
  &=& -\epsilon_{N-k}(\z_t)\mu_N(\z_t) \;=\; t\dot\zeta_{N-k}+o(t)
\end{eqnarray*}

Combining these equations with~(\ref{eq:Fzt}) and~(\ref{eq:tFzt}), we
find
\begin{eqnarray}\label{eq:phitphi}
\lefteqn{ -t^2F_1(\e)\Biggl\{\sum_{l=1}^N\dot\zeta_l\dot\zeta_{N-l}\biggl[1
  +{N\over 2}\sum_{d=1}^N\alpha^\2_d\omega_l^d\biggr]
  -N\alpha\dot\zeta_N^2\Biggr\}}\qqqq\nonumber\\
\lefteqn{\mbox{}\mbox{}+t^2N\sum_{k=1}^{N-1}\dot\zeta_k\dot\zeta_{N-k}
    \omega_k^{-m}\d{z_m}{\phi_k}(\e)
  +t^2N(1-N\alpha)\dot\zeta_N^2\d{z_m}{\phi_N}(\e)}\qquad\nonumber\\
&=& -t^2\wt F_1(\e)\Biggl\{\sum_{l=1}^N\dot\zeta_l\dot\zeta_{N-l}\biggl[1
  +{N\over 2}\sum_{d=1}^N[\alpha^\2_d-\t\alpha^\2_d]\omega_l^d\biggr]
  -N\alpha\dot\zeta_N^2\Biggr\}\nonumber\\
  & &\mbox{}+t^2N\sum_{k=1}^{N-1}\dot\zeta_k\dot\zeta_{N-k}
    \omega_k^{-m}\d{z_m}{\t\phi_k}(\e)\nonumber\\
& &\mbox{}+t^2N(1-N\alpha+N\t\alpha)\dot\zeta_N^2\d{z_m}{\t\phi_N}(\e)+o(t^2).
\end{eqnarray}

This equation in turn gives for $(1-N\alpha)\dot\zeta_N^2=1$ and
$\dot\zeta_l\dot\zeta_{N-l}=1$ if $1\leq l<N$
\begin{eqnarray*}
N\sum_{k=1}^N\omega_k^{-m}\d{z_m}{\phi_k}(\e)
&=& {N\t\alpha\over 1-N\alpha}N\d{z_m}{\t\phi_N}(\e)
    + N(F_1(\e)-\wt F_1(\e))\\
& &\mbox{}+{1\over2}\sum_{l=1}^{N-1}\sum_{d=1}^N\bar\alpha^\2_d\omega_l^d
          +{1\over1-N\alpha}{1\over2}\sum_{d=1}^N\bar\alpha^\2_d
\end{eqnarray*}

This expression, together with~(\ref{eq:dphiN})
and~(\ref{eq:mean}), yields Equation~(\ref{eq:dFe}). The derivation
of~(\ref{eq:sumdFe}) is done in a similar way: we use the fact that
\[
 \EE X_m 
  \;=\; \sum_{i=1}^N\d{z_m}{F_i}(\e)
  \;=\; N\d{z_m}{\phi_N}(\e).
\]

We use Equation~(\ref{eq:phitphi}) with $(1-N\alpha)\dot\zeta_N^2=1$ and
$\dot\zeta_k\dot\zeta_{N-k}=0$ for $k<N$ and find
\[
N\d{z_m}{\phi_N}(\e)
 = {1-N\alpha+N\t\alpha\over 1-N\alpha}N\d{z_m}{\t\phi_N}(\e)
    +F_1(\e)-\wt F_1(\e)
   +{1\over 1-N\alpha}{1\over2}\sum_{d=1}^N\bar\alpha^\2_d.
\]

This gives Equation~(\ref{eq:sumdFe}) and concludes the proof of the theorem.
\end{proof}

\begin{proof}{ of Theorem~\ref{thm:EW}}
We see easily that the mean numbers of clients arriving between
polling times are respectively $\alpha = \lambda\tau$ and $\t\alpha
= \lambda\t\tau$. Moreover, using Equations~(\ref{eq:al2})
and~(\ref{eq:tal2}) and the properties of generating functions we find
that, for $d<N$
\begin{eqnarray*}
\alpha_d^\2 &=& \alpha^\2 \;=\; \lambda^2\tau^\2\\
\alpha_N^\2 &=& \alpha^\2+\hat\lambda(b^\2-b)\tau\\
\t\alpha_d^\2 &=& \t\alpha^\2 \;=\; \lambda^2\t\tau^\2\\
\t\alpha_N^\2 &=& \t\alpha^\2+\hat\lambda(b^\2-b)\t\tau.
\end{eqnarray*}

The computation of waiting times uses the following classical
argument: a non empty queue visited by the server can be decomposed
into
\begin{itemize}
\item the head of line customer;
\item the clients who arrived after him during his waiting time;
\item clients who arrived in the same batch as the first client,
but are not yet served; since, by renewal arguments, the mean size of
this batch is $b^\2/b$, the mean number of remaining clients is
$(b^\2-b)/2b$.
\end{itemize}

This can be written as
\[
 \EE[X_m\mid S=m,X_m>0]=1+\lambda\EE[W]+{b^\2-b\over 2b}
\]
and, using the fact that
$\sum_{d=1}^N\alpha_d^\2\omega_l^d=\alpha_N^\2-\alpha^\2$,
\begin{eqnarray*}
 \lambda b\EE[W]
&=& {\EE[X_m\mid S=m]\over 1-P(X_m=0\mid S=m)}-1-{b^\2-b\over 2b}\\
&=& {\t\alpha\over 1-N\alpha}
       \sum_{l=1}^{N-1}{1\over 1-\mu_l}
    +{1\over 1-N\alpha}{N\over 2}\alpha^\2
    +{1\over N\t\alpha}{N\over 2}\t\alpha^\2\\
& &\mbox{}+{b^\2-b\over 2b}\Biggl[ 
            {\alpha+(N-1)\t\alpha\over1-N\alpha}
           -{N\t\alpha\over 1-N\alpha}{1\over N}
      \sum_{l=1}^{N-1}{\mu_l-\t\mu_l\over 1-\t\mu_l}\Biggr]\\
&=& {\lambda\t\tau\over 1-N\lambda\tau}
      \sum_{l=1}^{N-1}{1\over 1-\mu_l}
    +{N\lambda^2\tau^\2\over 2(1-N\lambda\tau)}
      +{\lambda\t\tau^\2\over2\t\tau}\\
& &\mbox{}+ {b^\2-b\over 2b}\Biggl[
      {\lambda\tau+(N-1)\lambda\t\tau\over 1-N\lambda\tau}
     -{\lambda\t\tau\over 1-N\lambda\tau}
\sum_{l=1}^{N-1}{\mu_l-\t\mu_l\over 1-\t\mu_l}\Biggr].
\end{eqnarray*}

This gives Equation~(\ref{eq:EW}).
\end{proof}

\bibliography{queues}
\bibliographystyle{siam}

\end{document}